\documentclass{amsart}

\theoremstyle{definition}

\theoremstyle{remark}

\usepackage{graphicx}
\usepackage{color}
\usepackage{amsmath}
\usepackage{amsfonts}
\usepackage{amssymb}
\newcommand{\be}{\begin{equation}}
\newcommand{\ee}{\end{equation}}

\newcommand{\lc}{{\stackrel{\scriptscriptstyle{LC}}{\Gamma}}}

\newcommand{\dz}{\wedge}

\newcommand{\ba}{\begin{array}}

\newcommand{\ea}{\end{array}}

\newcommand{\beq}{\begin{eqnarray}}

\newcommand{\eeq}{\end{eqnarray}}

\newtheorem{lm}{lemma}

\newtheorem{thee}{theorem}

\newtheorem{proo}{proposition}

\newtheorem{co}{corollary}

\newtheorem{rem}{remark}

\newtheorem{deff}{definition}

\newcommand{\bd}{\begin{deff}}

\newcommand{\ed}{\end{deff}}

\newcommand{\bl}{\begin{lm}}

\newcommand{\el}{\end{lm}}

\newcommand{\bp}{\begin{proo}}

\newcommand{\ep}{\end{proo}}

\newcommand{\bt}{\begin{thee}}

\newcommand{\et}{\end{thee}}

\newcommand{\bc}{\begin{co}}

\newcommand{\ec}{\end{co}}

\newcommand{\brm}{\begin{rem}}

\newcommand{\erm}{\end{rem}}

\newcommand{\der}{{\rm d}}

\newcommand{\sgn}{\mathrm{sgn}}

\usepackage{t1enc}
\def\frak{\mathfrak}

\newcommand{\newc}{\newcommand}

\let\ccdot\cdot
\def\cdot{\hbox to 2.5pt{\hss$\ccdot$\hss}}

\newc{\aR}{\mbox{\boldmath{$ R$}}}
\newc{\aS}{\mbox{\boldmath{$ S$}}}
\newc{\aT}{\mbox{\boldmath{$ T$}}}
\newc{\aW}{\mbox{\boldmath{$ W$}}}

\newc{\aK}{\mbox{\boldmath{$ K$}}}
\newc{\aL}{\mbox{\boldmath{$ L$}}}
\newcommand{\bbC}{\mathbb{C}}
\newcommand{\bbH}{\mathbb{H}}
\newcommand{\bbO}{\mathbb{O}}

\usepackage{amssymb}
\usepackage{amscd}

\newc{\obstrn}[2]{B^{#1}_{#2}}

\newcommand{\rpl}                         % +) or <+
{\mbox{$
\begin{picture}(12.7,8)(-.5,-1)
\put(0,0.2){$+$}
\put(4.2,2.8){\oval(8,8)[r]}
\end{picture}$}}

\newcommand{\lpl}                         % (+ or +>
{\mbox{$
\begin{picture}(12.7,8)(-.5,-1)
\put(2,0.2){$+$}
\put(6.2,2.8){\oval(8,8)[l]}
\end{picture}$}}

\usepackage{ifthen}

\newcommand{\bbR}{\mathbb{R}}
\newcommand{\sog}{\mathbf{SO}}

\newcommand{\og}{\mathbf{O}}
\newcommand{\soa}{\frak{so}}

\newcommand{\sua}{\frak{su}}
\newcommand{\spa}{\frak{sp}}

\newcommand{\spg}{\mathbf{Sp}}
\newcommand{\sug}{\mathbf{SU}}
\newcommand{\ug}{\mathbf{U}}

\newc{\tensor}[1]{#1}
\newc{\Mvariable}[1]{\mbox{#1}}
\newc{\down}[1]{{}_{#1}}
\newc{\up}[1]{{}^{#1}}

\newc{\JulyStrut}{\rule{0mm}{6mm}}
\newc{\midtenPan}{\mbox{\sf S}}
\newc{\midten}{\mbox{\sf T}}
\newc{\midtenEi}{\mbox{\sf U}}
\newc{\ATen}{\mbox{\sf E}}
\newc{\BTen}{\mbox{\sf F}}
\newc{\CTen}{\mbox{\sf G}}

\def\sideremark#1{\ifvmode\leavevmode\fi\vadjust{\vbox to0pt{\vss% the remark
 \hbox to 0pt{\hskip\hsize\hskip1em%                          will appear only
 \vbox{\hsize3cm\tiny\raggedright\pretolerance10000%          on the side
 \noindent #1\hfill}\hss}\vbox to8pt{\vfil}\vss}}}%
\newcommand{\bgw}{{\textstyle \bigwedge}}

\numberwithin{equation}{section}

\newcounter{romenumi}
\newcommand{\labelromenumi}{(\roman{romenumi})}

\newcommand{\ten}{\Upsilon}

\begin{document}
\title{Special Riemannian geometries modeled on the distinguished
symmetric spaces}

\author{Pawe\l~ Nurowski} \address{Instytut Fizyki Teoretycznej,
Uniwersytet Warszawski, ul. Hoza 69, Warszawa, Poland}
\email{nurowski@fuw.edu.pl} \thanks{This research was supported by
the KBN grant 1 P03B 07529}

\date{\today}

\begin{abstract} We propose studies of special Riemannian
geometries with structure groups $H_1=\sog(3)\subset\sog(5)$,
$H_2=\sug(3)\subset\sog(8)$, $H_3=\spg(3)\subset\sog(14)$ and
$H_4={\bf F_4}\subset\sog(26)$ in respective dimensions 5, 8, 14 and
26. These geometries, have torsionless models with symmetry groups
$G_1=\sug(3)$, $G_2=\sug(3)\times\sug(3)$, $G_3=\sug(6)$ and
$G_4={\bf E}_6$. The groups $H_k$ and $G_k$ constitute a part of the
`magic square' for Lie groups. Apart from the $H_k$ geometries in
dimensions $n_k$, the `magic square' Lie groups suggest studies of a
finite number of other special Riemannian geometries. Among them the
smallest dimensional are $\ug(3)$ geometries in dimension 12. The
other structure groups for these Riemannian geometries are: ${\bf
S}(\ug(3)\times\ug(3))$, $\ug(6)$, ${\bf E}_6\times\sog(2)$,
$\spg(3)\times\sug(2)$, $\sug(6)\times\sug(2)$,
$\sog(12)\times\sug(2)$ and ${\bf E}_7\times\sug(2)$. The respective
dimensions are: 18, 30, 54, 28, 40, 64 and 112. This list is
supplemented by the two `exceptional' cases of
$\sug(2)\times\sug(2)$ geometries in dimension 8 and
$\sog(10)\times\sog(2)$ geometries in dimension 32.
%Some homogeneous examples of such geometries in dimension 8 and
%... are given.

\vskip5pt\centerline{\small\textbf{MSC classification}: 53A40, 53B15,
  53C10}\vskip15pt
\end{abstract}
\maketitle
%*************
\section{Motivation}
The motivation for this paper comes from type II B string theory
(see e.g. \cite{nag}), where one considers $n=6$-dimensional compact
Riemannian manifold $(X,g)$ which, in addition to the Levi-Civita
connection $\nabla^{LC}$, is equipped with:
\begin{itemize}
    \item a {\it metric} conection $\nabla^T$ with {\it totally
    skew-symmetric torsion} $T$,
    \item a spinor field $\Psi$ on $X$.
\end{itemize}
{\it Special} Riemannian structure $(X,g,\nabla^T,T,\Psi)$ is
supposed to satisfy a number of field equations including:
$$
\nabla ^T\Psi=0,\quad \delta(T)=0,\quad T\cdot\Psi=\mu\Psi,\quad
Ric^{\nabla^T}=0.$$ To construct the solutions for these equations
one may proceed as follows.
\begin{itemize}
\item Let $\ten$ be an object (e.g. a tensor), whose isotropy under
the action of $\sog(n)$ is $H\subset \sog(n)$. Infinitesimaly, such
an object determines the inclusion of the Lie algebra $\mathfrak{h}$
of $H$ in $\soa(n)$.
\item If $(X,g)$ is endowed with such a $\ten$ we can decompose the
Levi-Civita connection 1-form $\lc\in\soa(n)\otimes\bbR^n$ onto
$\Gamma\in\frak{h}\otimes\bbR^n$ and the rest:
$$
\lc=\Gamma+\tfrac{1}{2}T.
$$
This decomposition is of course not unique, but let us, for a
moment, assume that we have a way of choosing it.
\item Then the first Cartan structure equation
$\der\theta+(\Gamma+\tfrac{1}{2}T)\dz\theta=0$ for the Levi-Civita
connection $\lc$ may be rewritten to the form
   $$\der\theta+\Gamma\dz\theta=-\tfrac12T\dz\theta$$
and may be interpreted as the first structure equation for a
   {\it metric} connection $\Gamma\in\mathfrak{h}\otimes\bbR^n$ with
   torsion $T\in\soa(n)\otimes\bbR^n$.
  \item Curvature of this connection $K\in\mathfrak{h}\otimes\bigwedge^2\bbR^n$
  is determined via the second structure equation:
   $$K=\der\Gamma+\Gamma\dz\Gamma.$$
\item To escape from the ambiguity in the split $\lc=\Gamma+\tfrac12
T$ we impose the type II string theory requirement, that $T$ should
be totally skew-symmetric. \item Thus we require that our special
Riemannian geometry $(X,g,\ten)$ must admit a split
$\lc=\Gamma+\tfrac{1}{2}T$ with $T\in \bigwedge^3\bbR^n$ and
$\Gamma\in\mathfrak{h}\otimes\bbR^n$.
\item Examples are known of special Riemannian geometries (e.g.
nearly K\"ahler geometries in dimension $n=6$) in which such a
requirement uniquely determines $\Gamma$ and $T$.
\end{itemize}

This leads to the following problem: find all geometries
$(X,g,\ten)$ admitting the {\it unique} split $$
\lc=\Gamma+\tfrac{1}{2}T\quad {\rm with}\quad T\in
\bigwedge^3\bbR^n\quad{\rm and}\quad
\Gamma\in\mathfrak{h}\otimes\bbR^n. $$ In which dimensions $n$ do they
exist? What is $\ten$ which reduces $\sog(n)$ to $H$ for them? What
are the possible isotropy groups $H\in \sog(n)$?

\section{Special geometries $(X,g,\nabla^T,T,\Psi)$ }
If $T\in\bgw^3\bbR^n$ was {\it identically zero}, then since
$\mathfrak{h}\otimes\bbR^n\ni\Gamma=\lc$, the {\it holonomy group}
of $(X,g)$ would be {\it reduced} to $H\subset\sog(n)$. {\it All
irreducible} compact Riemannian manifolds $(X,g)$ with the reduced
holonomy group are classified by Berger \cite{ber}. These are:
\begin{itemize}\item either {\it symmetric spaces} $G/H$, with the
holonomy group $H\subset \sog(n)$\item or they are contained in the
{\it Berger's list}:
\end{itemize}
\begin{tabular}{|l|l|l|l|}
\hline Holonomy group for $g$ &Dimension of $X$&Type of $X$&Remarks\\
\hline
$\sog(n)$&$n$&generic&\\
&&&\\
$\ug(n)$&$2n$, $n\geq 2$&K\"ahler manifold&K\"ahler\\
&&&\\
$\sug(n)$&$2n$, $n\geq 2$&Calabi-Yau manifold&Ricci-flat,K\"ahler \\
&&&\\
$\spg(n)\cdot \spg(1)$&$4n$, $n\geq 2$&quaternionic K\"ahler&Einstein\\
&&&\\
$\spg(n)$&$4n$, $n\geq 2$&hyperk\"ahler manifold&Ricci-flat,K\"ahler\\
&&&\\
${\bf G}_2$&$7$&${\bf G}_2$ manifold&Ricci-flat\\
&&&\\
${\bf Spin}(7)$&$8$&${\bf Spin}(7)$ manifold&Ricci-flat\\
%&&&\\
%$Spin(9)$&$16$&$Spin(9)$ manifold&Einstein\\
% wyeliminowane przez Aleksejewskiego Graya i Brawna (szarego i brazowego)
\hline
\end{tabular}

\vspace{1cm}
We may relax $T=0$ for geometries with $H$ from Berger's theorem in at
least two ways:
\begin{itemize}
\item relax $T=0$ condition to $T\in\bigwedge^3\bbR^n$ for $H$ from {\it the Berger's
list}. This approach leads e.g. to {\it nearly} K\"ahler geometries
for $H=U(n)$, special {\it nonintegrable} $SU(3)$ geometries in
dimension $6$, special {\it nonintegrable} $G_2$ geometries in
dimension $7$, etc.
\item relax $T=0$ condition to $T\in\bigwedge^3\bbR^n$ for $H$
corresponding to the irreducible symmetric spaces $G/H$ from {\it
Cartan's list}.
\end{itemize}
In this note we focus on the second possibility. Here the simplest
case is related to the first entry in Cartan's list of the
irreducible symmetric spaces, namely to $G/H=\sug(3)/\sog(3).$ Thus, one
considers a $n=5$-dimensional manifold $X=\sug(3)/\sog(3)$ with the
{\it irreducible} $\sog(3)$ action at each tangent space at every
point of $X$. In such an approach $X=\sug(3)/\sog(3)$ is the {\it
integrable} ($T=0$) model for the irreducible $\sog(3)$ geometries
in dimension 5.

\vspace{0.5cm}
 Friedrich \cite{fried} asked the following questions: is it
possible to have $5$-dimensional Riemannian geometries
$(X,g,\nabla^T,T,\Psi)$ for which the torsionless model would be
$G/H=\sug(3)/\sog(3)$? If so, what is $\ten$ for such geometries?
\vspace{0.5cm}

\noindent In a joint work with Bobie\'nski \cite{bobi} we answered 
these questions as follows: \vspace{0.2cm}
\begin{itemize}
\item Tensor $\ten$ whose isotropy group under the action of $\sog(5)$ is the irreducible
$\sog(3)$ is determined by the following conditions:
\begin{itemize}
\item[i)] %it was totally symmetric:
$\ten_{ijk}=\ten_{(ijk)},$\hspace{2cm}(totally {\it symmetric})
\item[ii)]%it was trace free:
$\ten_{ijj}=0,$\hspace{3.5cm}(trace-free)
\item[iii)] $\ten_{jki}\ten_{lmi}+\ten_{lji}\ten_{kmi}+\ten_{kli}\ten_{jmi}=g_{jk}g_{lm}+g_{lj}g_{km}+g_{kl}g_{jm}.$
\end{itemize}
\item A $5$-dimensional Riemannian manifold $(X,g)$
equipped with a tensor field $\ten$ satisfying conditions i)-iii)
and admitting a unique decomposition $\lc=\Gamma+\tfrac{1}{2}T,$
with $T\in\bigwedge^3\bbR^5$ and $\Gamma\in\soa(3)\otimes\bbR^5$ is
called {\it nearly integrable} irreducible $\sog(3)$ structure.
\item We have examples of such geometries. All our examples admit transitive
symmetry group (which may be of dimension $8$, $6$ and $5$)
\item In particular, we have a $7$-parameter family of nonequivalent
examples which satisfy
$$
\nabla ^T\Psi=0,\quad\delta(T)=0,\quad T\cdot\Psi=\mu\Psi
$$
i.e. equations of type IIB string theory (but in wrong dimension!).
For this family of examples $T\neq 0$ and, at every point of $X$, we
have two $2$-dimensional vector spaces of $\nabla^T$-covariantly
constant spinors $\Psi$. Moreover, since for this family $K=0$, we
also have $Ric^{\nabla^T}=0$.
\end{itemize}

\section{Distinguished dimensions}
A natural question is \cite{nur}: what are the possible dimensions
$n$ in which there exists a tensor $\ten$ satisfying:
\begin{itemize}
\item[i)] %it was totally symmetric:
$\ten_{ijk}=\ten_{(ijk)},$\hspace{2cm} (total {\it symmetry})
\item[ii)]%it was trace free:
$\ten_{ijj}=0,$\hspace{3.5cm}({\it no} trace)
\item[iii)] $\ten_{jki}\ten_{lmi}+\ten_{lji}\ten_{kmi}+\ten_{kli}\ten_{jmi}=
g_{jk}g_{lm}+g_{lj}g_{km}+g_{kl}g_{jm}$?
\end{itemize}\vspace{0.2cm}
In dimension $n=5$ tensor $\ten$ has the following features
\cite{bobi}:\vspace{0.2cm}

\begin{itemize}
\item Given $\ten_{ijk}$ we consider a $3$rd order polynomial
 $w(a)=\ten_{ijk}a_ia_ja_k$, where
$a_i\in\bbR$, $i=1,2,3,4,5$. We have:
$$w(a)=6\sqrt{3}a_1a_2a_3+3\sqrt{3}(a_1^2-a_2^2)a_4-\big(~3a_1^2+3a_2^2-6a_3^2-6a_4^2+2a_5^2~\big)a_5
$$
\item Note that:
\be w(a)=\det
\begin{pmatrix}
a_5-\sqrt{3}a_4&\sqrt{3}a_3&\sqrt{3}a_2\\
\sqrt{3}a_3&a_5+\sqrt{3}a_4&\sqrt{3}a_1\\
\sqrt{3}a_2&\sqrt{3}a_1&-2a_5
\end{pmatrix}\label{det}
\ee
\end{itemize}
The last observation led Bryant \cite{bry} to the following answer
to our question from the beginning of this section: if tensor $\ten$
with the properties i)-iii) exists in dimension $n=5$, then it also
exists in dimensions $n=8$, $n=14$ and $n=26$. This is because in
addition to the field of the real numbers $\bbR$, we also have
$\bbC$, $\bbH$ and $\bbO$.

We now change expression (\ref{det}) into
\be
w(a)=\det A=\det\begin{pmatrix}
a_5-\sqrt{3}a_4&\sqrt{3}{\bf \alpha_3}&\sqrt{3}{\bf \alpha_2}\\
\sqrt{3}{\bf \overline{\alpha}_3}&a_5+\sqrt{3}a_4&\sqrt{3}{\bf \alpha_1}\\
\sqrt{3}{\bf \overline{\alpha}_2}&\sqrt{3}{\bf
\overline{\alpha}_1}&-2a_5
\end{pmatrix}\label{wa}\ee
where $a_i\in\bbR$, $i=1,2,3,...n$ and\begin{itemize}
\item for
$n=5$ we have:
\begin{eqnarray*}
{\bf\alpha_1}&=&a_1 \\
{\bf\alpha_2}&=&a_2 \\
{\bf\alpha_3}&=&a_3
%{\color{white}dla} $n=5$\\
\end{eqnarray*}
\item for $n=8$ we have
\begin{eqnarray*}
{\bf\alpha_1}&=&a_1+a_6 {\rm i} \\
{\bf\alpha_2}&=&a_2+a_7 {\rm i} \\
{\bf\alpha_3}&=&a_3+a_8 {\rm i}
%{\color{white}dla} $n=5$\\
\end{eqnarray*}
\item for $n=14$ we have:
\begin{eqnarray*}
{\bf\alpha_1}&=&a_1+a_6 {\rm i}+a_9 {\rm j}+a_{10}{\rm k} \\
{\bf\alpha_2}&=&a_2+a_7 {\rm i}+a_{11} {\rm j}+a_{12}{\rm k} \\
{\bf\alpha_3}&=&a_3+a_8 {\rm i}+a_{13} {\rm j}+a_{14}{\rm k}
%{\color{white}dla} $n=5$\\
\end{eqnarray*}
\item for $n=26$ we have:
\begin{eqnarray*}
{\bf\alpha_1}&=&a_1+a_6 {\rm i}+a_9 {\rm
j}+a_{10}{\rm k} +a_{15} {\rm p}+a_{16} {\rm q}+a_{17}
{\rm r}+a_{18}{\rm s}\\
{\bf\alpha_2}&=&a_2+a_7 {\rm i}+a_{11} {\rm j}+a_{12}{\rm k}+a_{19}
{\rm p}+a_{20} {\rm q}+a_{21}
{\rm r}+a_{22}{\rm s} \\
{\bf\alpha_3}&=&a_3+a_8 {\rm i}+a_{13} {\rm j}+a_{14}{\rm k}+a_{23}
{\rm p}+a_{24} {\rm q}+a_{25} {\rm r}+a_{26}{\rm s}
%{\color{white}dla} $n=5$\\
\end{eqnarray*}
\end{itemize}
with ${\rm i}$ imaginary unit, ${\rm i,j,k}$ imaginary quaternion
units and ${\rm i,j,k,p,q,r,s}$ imaginary octonion units.
Remarkably, modulo the action of the $\og(n)$ group, the symbol
$\det$ in the above expression is well defined by the demand that the 
Weierstrass formula
$$\det A=\sum_{\pi\in S_3} \sgn\pi~
A_{1\pi(1)}A_{2\pi(2)}A_{3\pi(3)}$$ assumes {\it real} values.

Now, we have the following theorems \cite{car1,car2}.\\
{\bf Theorem 1}\\
For each $n=5,8,14$ i $26$ tensor $\ten$ given by
$$\ten_{ijk}a_ia_ja_k=w(a)=\det A$$
 satisfies i)-iii).\\
{\bf Theorem 2}\\
In dimensions $n=5$, $8$, $14$ i $26$ tensor $\ten$ reduces the
${\bf GL}(n,\bbR)$ group via $\og(n)$ to a subgroup
  $H_n$, where:
  \begin{itemize}
  \item for $n=5$ ~~group $H_5$ ~is the irreducible $\sog(3)$ in
  $\sog(5)$;\\ Here, the torsionless compact model is: $\sug(3)/\sog(3)$
  \item for $n=8$ ~~group $H_8$ ~is the irreducible $\sug(3)$ in
  $\sog(8)$;\\Here, the torsionless compact model is: $\sug(3)$
  \item for $n=14$ group $H_{14}$ is the irreducible $\spg(3)$ ~in
  $\sog(14)$;\\ Here, the torsionless model is: $\sug(6)/\spg(3)$
  \item for $n=26$ group $H_{26}$ is the irreducible ${\bf F}_4$ ~~~~~~in $\sog(26)$;
  \\ Here, the torsionless compact model is:
  ${\bf E}_6/{\bf F}_4$
  \end{itemize}
{\bf Theorem 3}
\begin{itemize}
\item The only dimensions in which conditions i)-iii) have solutions for $\ten_{ijk}$ are
$n=5,8,14,26$.
\item Modulo the action of $\og(n)$ all such tensors
are given by $\det A$, where $A$ is a $3\times 3$ traceless
hermitian matrix with entries in $\bbR,\bbC,\bbH,\bbO$, for the
respective dimensions $5$,$8$,$14$,$26$.
\end{itemize}
\vspace{0.5cm}

\noindent {\bf Idea of the proofs:}
\begin{itemize}
\item The theorems follow from Cartan's works on {\it isoparametric hypersurfaces in
spheres} \cite{car1,car2}.
\item A hypersurface $S$ is isoparametric in ${\bf S}^{n-1}$ iff all its
{\it principal curvatures} are {\it constant}.
\item Cartan proved that $S$ is isoparametric in
$${\bf S}^{n-1}=\{a^i\in \bbR^n~|~
(a^1)^2+(a^2)^2+...+(a^n)^2=1\}$$ and has $3$ {\it distinct}
principal curvatures iff $S={\bf S}^{n-1}\cap P_c$,
where
$$P_c=\{a^i\in \bbR^n~|~w(a)=c=const\in\bbR\}$$ and $w=w(a)$ is
a homogeneous $3$rd order {\it polynomial} in variables $(a^i)$ such
that
\begin{eqnarray*} {\rm cii)}&&\triangle w=0\\
{\rm ciii)}&&|\nabla w|^2=9~[~(a^1)^2+(a^2)^2+...+(a^n)^2~]^2.
\end{eqnarray*}
\item He reduced the above differential equations for $w=w(a)$ to
equations for a certain function with the properties of a function
he encountered when solving the problem of finding Riemannian spaces
with absolute parallelism. He proved that such function give rise only
  to $w=w(a)$ given by formula (\ref{wa}). Of course, Cartan's
  conditions cii)-ciii) for the polynomial $w=w(a)$ translate to our
  conditions ii)-iii) for the corresponding symmetric tensor $\ten$
  such that $w(a)=\ten_{ijk}a^ia^ja^k$.
\end{itemize}
These three theorems lead to the idea \cite{nur} of studies of $H_k$
structures
in dimensions $n_k=5,8,14,26$.\\
{\bf Definition}\\An $H_k$ structure on a $n_k$-dimensional
Riemannian manifold $(M,g)$ is a structure defined by means of a
rank 3 tensor field $\ten$ satisfying
\begin{enumerate}
\item[i)] $\ten_{ijk}=\ten_{(ijk)}$,
\item[ii)]%it was trace free:
$\ten_{ijj}=0$, \vskip1mm
\item[iii)]
  $\ten_{jki}\ten_{lmi}+\ten_{lji}\ten_{kmi}+\ten_{kli}\ten_{jmi}=g_{jk}g_{lm}
  +g_{lj}g_{km}+g_{kl}g_{jm}.$
\end{enumerate}\vspace{0.1cm}
An $H_k$ structure is called {\it nearly integrable} iff
$$\nabla^{LC}_X\ten(X,X,X)=0,\quad\quad\forall X\in\Gamma(TM).$$

\section{Nearly integrable $H_k$ structures and characteristic
connection} Now a natural question is: what are the necessary and
sufficient conditions for a $H_k$ structure to admit a unique
decomposition $$\lc=\Gamma+\tfrac12 T\quad{\rm with}\quad
\Gamma\in\mathfrak{h}_k\otimes\bbR^k\quad{\rm and}\quad
T\in\bgw^3\bbR^{n_k}?$$

If such a unique decomposition exists, the connection $\Gamma$ is
called {\it characteristic connection} of the $H_k$ structure.\\

The answer to the above question is given by the following theorem \cite{nur}.\\
{\bf Theorem 4}\\
Every $H_k$ structure that admits a characteristic connection must
be {\it nearly integrable}.\\Moreover,
\begin{itemize} \item In dimensions $5$ and $14$ the nearly
integrable condition is also sufficient for the existence of the
characteristic connection.\item In dimension $8$ the spaces
$\mathfrak{h}_k\otimes\bbR^k$ and $\bgw^3\bbR^{n_k}$ have
$1$-dimensional intersection $V_1$. In this dimension a sufficient
condition for the existence of characteristic connection $\Gamma$ is
that the Levi-Civita connection $\lc$ of a nearly integrable
$\sug(3)$ structure does not have $V_1$ components in the $\sug(3)$
decomposition of $\soa(8)\otimes\bbR^8$ onto the irreducibles.\item
In dimension $26$ the Levi-Civita connection $\lc$ of a nearly
integrable ${\bf F}_4$ structure may have values in $52$-dimensional
irreducible representation $V_{52}$ of ${\bf F}_4$, which is not
present in the algebraic sum of $\mathfrak{f}_4\otimes\bbR^k$ and
$\bgw^3\bbR^{n_k}$. The sufficient condition for such structures to
admit characteristic $\Gamma$ is that $\lc$ has not components in
$V_{52}$.
\end{itemize}
{\bf Definition}\\The nearly integrable $H_k$ structures described
by Theorem 4 are called {\it restricted} nearly integrable.

Now we discuss what the restricted nearly integrable condition means
for a $H_k$ structure:
\begin{itemize}
\item If $n_k=5$ then, out of the {\it a priori} $50$ independent
components of the Levi-Civita connection $\lc$, the restricted
nearly integrable condition excludes $25$. Thus, heuristically, the
restricted nearly integrable $\sog(3)$ structures constitute `a
half' of all the possible $\sog(3)$ structures in dimension $5$.
\item If $n_k=8$ the Levi-Civita connection has $224$ components. The
restricted nearly integrable condition reduces it to $118$.
\item For
$n_k=14$ these numbers reduce from $1274$ to $658$.
\item For
$n_k=26$ the reduction is from $8450$ to $3952$.
\end{itemize}
To discuss the possible torsion types of the characteristic
connection for $H_k$ geometries we need to know that: \begin{itemize}
\item there are {\it real} irreducible
representations of the group $\sog(3)$ in odd dimensions:
$1,3,5,7,9...$
\item there are {\it real} irreducible
representations of the group $\sug(3)$ in dimensions:
$1,8,20,27,70...$ \item there are {\it real} irreducible
representations of the group $\spg(3)$ in dimensions:
$1,14,21,70,84,90,126,189,512,525...$
\item there are {\it
real} irreducible representations of the group ${\bf F}_4$ in
dimensions: $1,26,52,273,324,1053,1274,4096,8424...$
\end{itemize}
For each of the possible dimensions $n_k=5,8,14,26$ we denote a
possible irreducible $j$-dimensional representation of the $H_k$
group by ${^{n_k}V_j}$. Then we have the following theorem \cite{nur}.\\
{\bf Theorem 5}\\
 Let
$(M,g,\ten)$ be a nearly integrable $H_k$ structure admitting
characteristic connection $\Gamma$. The $H_k$ irreducible
decomposition of the skew symmetric torsion $T$ of $\Gamma$ is given
by:
\begin{itemize}
\item
$T\in {^5V_7}\oplus{^5V_3},$ \hspace{3.5cm}for $n_k=5$,
\item
$T\in
  {^8V_{27}}\oplus{^8V_{20}}\oplus{^8V_{8}}\oplus{^8V_{1}},$
  \hspace{1.37cm}for $n_k=8$,
\item $T\in
  {^{14}V_{189}}\oplus{^{14}V_{84}}\oplus{^{14}V_{70}}\oplus{^{14}V_{21}},$
 \hspace{.36cm}for $n_k=14$,
\item
$T\in {^{26}V_{1274}}\oplus{^{26}V_{1053}}\oplus{^{26}V_{273}},$
 \hspace{1.1cm}for $n_k=26$.
\end{itemize}
{\bf Example \cite{nur}: $\sug(3)$ structures in dimension $8$ }\\
Among many interesting features of these structures we mention the
following:
\begin{itemize}
\item We
have examples of these structures admitting a characteristic
connection with nonzero torsion.
\item All our examples admit transitive symmetry group, which can
has dimension $\leq 16$.
\item We have $2$-parameter family of examples
with transitive symmetry group of dimension $11$, with torsion
$T\in{^8V_{27}}$ and the Ricci tensor $Ric^\Gamma$ of the
characteristic connection $\Gamma$ with 2 different constant
eigenvalues of multiplicity $5$ and $3$ \item We have another
$2$-parameter family of examples with transitive symmetry group of
dimension $9$, with {\it vectorial} torsion $T\in{^8V_{8}}$ and with
the Ricci tensor $Ric^\Gamma$ of the characteristic connection
$\Gamma$ with 2 different constant eigenvalues of multiplicity $4$
and $4$.
\item
In the decomposition of $\bgw^3\bbR^8$ onto the irreducible
components under the action of $\sug(3)$ there exists a
$1$-dimensional $\sug(3)$ invariant subspace ${^8V}_1$.
\item This space, in an orthonormal coframe adapted to the $\sug(3)$
  structure, 
is spanned by a $3$-form $$\psi=\tau_1\dz\theta^6+\tau_2\dz\theta^7
+\tau_3\dz\theta^8+\theta^6\dz\theta^7\dz\theta^8,$$ where
$(\tau_1,\tau_2,\tau_3)$ are $2$-forms
\begin{eqnarray*}
\tau_1&=&\theta^1\dz\theta^4+\theta^2\dz\theta^3+\sqrt{3}\theta^1\dz\theta^5
\\
\tau_2&=&\theta^1\dz\theta^3+\theta^4\dz\theta^2+\sqrt{3}\theta^2\dz\theta^5\\
\tau_3&=&\theta^1\dz\theta^2+2\theta^4\dz\theta^3
\end{eqnarray*}
spanning the $3$-dimensional irreducible representation
${^5\bgw}^2_3\simeq\soa(3)$ associated with $\sog(3)$ structure in
dimension $5$.
\item The $3$-form $\psi$ can be considered in $\bbR^8$ without
any reference to tensor $\ten$.\item It is remarkable that this
$3$-form {\it alone} reduces the ${\bf GL}(8,\bbR)$ to the
irreducible $\sug(3)$ in the same way as $\ten$ does.
\item Thus, in
dimension $8$, the $H_k$ structure can be defined either in terms of
the {\it totally symmetric} $\ten$ or in terms of the {\it totally
skew symmetric} $\psi$. \item In this sense the 3-form $\psi$ and
the 2-forms $(\tau_1,\tau_2,\tau_3)$ play the same role in the
relations between $\sug(3)$ structures in dimension {\it eight} and
$\sog(3)$ structures in dimension {\it five} as the 3-form
$$\phi=\sigma_1\dz\theta^5+\sigma_2\dz\theta^6
+\sigma_3\dz\theta^7+\theta^5\dz\theta^6\dz\theta^7$$ and the
self-dual 2-forms
\begin{eqnarray*}
\sigma_1&=&\theta^1\dz\theta^3+\theta^4\dz\theta^2\\
\sigma_2&=&\theta^4\dz\theta^1+\theta^3\dz\theta^2\\
\sigma_3&=&\theta^1\dz\theta^2+\theta^3\dz\theta^4
\end{eqnarray*}
play in the relations between ${\bf G}_2$ structures in dimension
{\it seven} and $\sug(2)$ structures in dimension {\it four}.
\end{itemize}
\section{The magic square and special geometries modeled on distinguished
symmetric spaces} Looking at the torsionless models $X_k$ for the
$H_k$ structures given in Theorem 2 we see that each group $H_k$ has
its associated group $G_k$ such that $X_k=G_k/H_k$. Remarkably the
Lie algebras of the groups $H_k$ and $G_k$ constitute, respectively,
the first and the second column of the celebrated `magic square' of
the Lie algebras \cite{mag1,mag2}:
\begin{center}
\begin{tabular}{|c|c|c|c|}
\hline
$\soa(3)$&$\sua(3)$&$\spa(3)$&$\mathfrak{f}_4$\\
\hline
$\sua(3)$&$2\sua(3)$&$\sua(6)$&$\mathfrak{e}_6$\\
\hline
$\spa(3)$&$\sua(6)$&$\soa(12)$&$\mathfrak{e}_7$\\
\hline
$\mathfrak{f}_4$&$\mathfrak{e}_6$&$\mathfrak{e}_7$&$\mathfrak{e}_8$\\
\hline
\end{tabular}~~.
\end{center}
Let $G_k$, $\mathcal{G}_k$ and $\tilde{\mathcal{G}}_k$ denote the
respective compact Lie groups corresponding to the Lie algebras of
the second, third and the fourth columns of the magic
square. The observation opening this section suggests that the pairs
$(G_k,\mathcal{G}_k)$ and $(\mathcal{G}_k,\tilde{\mathcal{G}}_k)$,
with the homogeneous spaces $\mathcal{G}_k/G_k$ and
$\tilde{\mathcal{G}}_k/\mathcal{G}_k$ may model $T\equiv 0$ cases of
other interesting special Riemannian geometries with skew-symmetric
torsion. This suggestion should be a bit modified, since the
homogeneous spaces $\mathcal{G}_k/G_k$ and
$\tilde{\mathcal{G}}_k/\mathcal{G}_k$ are reducible. To have {\it
irreducible} symmetric spaces we need \begin{itemize} \item either
to take the Lie groups $\mathcal{G}_k$ corresponding to the third
column of the magic square and divide them by the compact Lie groups
corresponding to the Lie algebras from the following table:
\begin{center}
\begin{tabular}{|c|}
\hline
$\sua(3)\oplus\bbR$\\
\hline
$2\sua(3)\oplus\bbR$\\
\hline
$\sua(6)\oplus\bbR$\\
\hline
$\mathfrak{e}_6\oplus\bbR$\\
\hline
\end{tabular}~~,
\end{center}
 \item or we need to take the Lie groups $\tilde{\mathcal{G}_k}$
corresponding to the forth column of the magic square and divide
them by the compact Lie groups corresponding to the Lie algebras
from the following table:
\begin{center}
\begin{tabular}{|c|}
\hline
$\spa(3)\oplus\sua(2)$\\
\hline
$\sua(6)\oplus\sua(2)$\\
\hline
$\soa(12)\oplus\sua(2)$\\
\hline
$\mathfrak{e}_7\oplus\sua(2)$\\
\hline
\end{tabular}~~.
\end{center}
\end{itemize}
This leads to twelve torsionless models of special Riemannian
geometries \cite{nur} given in the table below:
\begin{center}
\begin{tabular}{|c|c|c|}
\hline
$\sug(3)/\sog(3)$&$\spg(3)/\ug(3)$&${\bf F}_4/(\spg(3)\times\sug(2))$\\
\hline
$\sug(3)$&$\sug(6)/{\bf S}(\ug(3)\times\ug(3))$&${\bf E}_6/(\sug(6)\times\sug(2)$\\
\hline
$\sug(6)/\spg(3)$&$\sog(12)/\ug(6)$&${\bf E}_7/(\sog(12)\times\sug(2))$\\
\hline
${\bf E}_6/{\bf F}_4$&${\bf E}_7/({\bf E}_6\times\sog(2))$&${\bf E}_8/({\bf E}_7\times\sug(2))$\\
\hline
\end{tabular}~~.
\end{center}
It is an interesting question if these 12 symmetric spaces can be
deformed to obtain twelve classes of special geometries $X$ with
totally skew symmetric torsion and the characteristic connection.
The dimensions $n$ of $X$ and the structure groups of the
characteristic connection for these geometries are given in the
table below:
\begin{center}
\begin{tabular}{||c|c||c|c||c|c||}
\hline
 $n=$&Structure&$n=$&Structure group&$n=$&Structure\\
 $n_k$& group $H_k$&$2(n_k+1)$&&$4(n_k+2)$&group\\
\hline\hline
5&$\sog(3)$&12&$\ug(3)$&28&$\spg(3)\times\sug(2)$\\\hline
8&$\sug(3)$&18&${\bf
S}(\ug(3)\times\ug(3))$&40&$\sug(6)\times\sug(2)$\\\hline
14&$\spg(3)$&30&$\ug(6)$&64&$\sog(12)\times\sug(2)$\\\hline 26&${\bf
F}_4$&54&${\bf E}_6\times\sog(2)$&112&${\bf E}_7
\times\sug(2)$\\\hline
\end{tabular}
\end{center}
A quick look at the Cartan's list of symmetric spaces shows that
this list should be supplemented by two exceptional cases:
\begin{itemize}
\item[1)] $\dim X=8$, with the structure
group $\sug(2)\times\sug(2)$ and with the torsionless model of
compact type $X={\bf G}_2/(\sug(2)\times\sug(2))$ .
\item[2)] $\dim X=32$, with the structure group $\sog(10)\times\sog(2)$
and with the torsionless model of compact type $X={\bf
E}_6/(\sog(10)\times\sog(2))$
\end{itemize}
Besides the geometries from the first column of the above table, and
besides the exceptional geometries of case 1) above (see \cite{nur}),
we do not know what objects $\ten$ reduce the $\og(n)$ groups to the
structure groups included in the table.
\section{Acknowledgements}
I wish to thank Tohru Morimoto for inviting me to the `RIMS
Symposium on developments of Cartan geometry and related
mathematical problems' and Ilka Agricola for inviting me to the
workshop `Special geometries in mathematical physics' held in
K\"uhlungsborn. The present work was initiated at the RIMS workshop
in October 2005, where Robert Bryant informed me
about relations between tensor $\ten$ of Ref. \cite{bobi} and
Cartan's works on isoparametric hypersurfaces. This resulted in the
paper \cite{nur}. The present work is a compacted version of
\cite{nur} prepared for a talk which I gave in K\"uhlungsborn in
March 2006.

\end{document}